\documentclass[preprint,review,10pt]{elsarticle}
\usepackage{amssymb}
\usepackage{pifont}
\usepackage{amsfonts}
\usepackage{amsmath,amssymb}
\usepackage{graphicx}
\usepackage{setspace}
\usepackage[left=2cm,top=1.75cm,right=2cm]{geometry}
\usepackage{subfig}
\usepackage{natbib}

\newtheorem{corollary}{Corollary}

\newtheorem{lemma}{Lemma}

\newtheorem{remark}{Remark}

\newtheorem{theorem}{Theorem}

\numberwithin{equation}{section} \journal{XYZ}
\begin{document}
\begin{frontmatter}
\title{Asymptotic formula for $q$-Derivative of $q$-Durrmeyer Operators}
\author[label a,label b,label *]{Prashantkumar Patel}
\ead{prashant225@gmail.com}
\author[label a,label d]{Vishnu Narayan Mishra}
\ead{vishnu\_narayanmishra@yahoo.co.in; vishnunarayanmishra@gmail.com}
\author[label c]{R. N. Mohapatra}
\ead{ramm@pegasus.cc.ucf.edu}
\address[label a]{Department of Applied Mathematics \& Humanities,
S. V. National Institute of Technology,
Surat-395 007 (Gujarat), India}
\address[label c]{Mathematics Department, University of Central Florida, Orlando, FL 32816, USA}
\address[label b]{Department of Mathematics, St. Xavier College, Ahmedabad-380 009 (Gujarat), India}
\address[label d]{ L. 1627 Awadh Puri Colony Beniganj, Phase -III, Opposite - Industrial Training Institute (I.T.I.), Ayodhya Main Road,
Faizabad-224 001, (Uttar Pradesh), India}
\fntext[label*]{Corresponding author}
\begin{abstract}
In the manuscript, Voronovskaja type asymptotic formula for function having $q$-derivative of $q$-Durrmeyer operators and
$q$-Durrmeyer-Stancu operators are discussed.
\end{abstract}
\begin{keyword} $q$-integers; $q$-Durrmeyer operators; $q$-derivative; asymptotic formula. \\
\textit{2000 Mathematics Subject Classification: } 41A25 · 41A28 · 41A35 · 41A36 \end{keyword}
\end{frontmatter}
\section{Introduction}
The classical Bernstein-Durrmeyer operators $D_n$ introduced by
Durrmeyer \cite{Durrmeyer1967} associate with each function $f$
integrable on the interval $[0, 1]$, the polynomial
\begin{equation}\label{12.eq1.1}
D_n(f;x) = (n+1) \sum_{k=0}^{n} p_{n,k} (x) \int_0^1 p_{n,k}(t) f(t) dt, ~~~ x\in [0,1],
\end{equation}
where $\displaystyle p_{n,k}(x) = {n \choose k} x^k(1-x)^{n-k}$. \\
\indent These operators been studied by Derriennic \cite{MMD} and many others. Last 30 years, the application of
$q$-calculus in filed of approximation theory is active area of research. In 1987,
the $q$-analogues of Bernstein operators was introduced by Lupas \cite{lupas1987q},
 Gupta and Hapeing \cite{gupta2008rate} introduced $q$-generalization of the operators \eqref{12.eq1.1} as
\begin{equation}\label{12.eq1.2}
D_{n,q}(f; x)=[n+1]_q \sum_{k=0}^n q^{-k} p_{nk}(q; x)\int_0^1 f(t)p_{nk}(q; qt)d_qt,
\end{equation}
where  $\displaystyle p_{nk}(q;x) = {n\choose k}_q x^k(1-x)_q^{n-k}$.\\
\indent The Rate of convergence of the operators \eqref{12.eq1.2}
was discussed by Gupta \textit{et al.}
\cite{Gupta2008some,zeng2010note}, local approximation, global
approximation and simultaneous approximation properties of these
operators by Finta and Gupta \cite{finta2009approximation},
estimation  of moments and King type approximation was elaborated
by Gupta and Sharma \cite{gupta2011recurrence}. In 2014, Mishra
and Patel \cite{mishra2013short,mishra2014generalized} talk about
Stancu generalization, Voronovskaja type asymptotic formula and
various other approximation properties of the $q$-Durrmeyer-Stancu
operators. We have the notation of $q$-calculus as given in
\cite{KacCheung2002,thomae1869beitrage}. Here, in this manuscript
we establish  Voronovskaja type asymptotic formula for function
having $q$-derivative.
\section{Estimation of moments and Asymptotic formula}
In the sequel, we shall need the following auxiliary results:
\begin{theorem}\cite{gupta2011recurrence}
If $m$-th $(m >0,m\in \mathbb{N})$ order moments of operator \eqref{12.eq1.2} is defined as
$$D^q_{n,m}(x) = D_{n,q}(t^m, x) = [n + 1]_q \sum_{k=0}^n q^{-k}p_{n,k}(q; x)\int_0^1 p_{n,k}(q; qt)t^md_qt, ~~~x\in [0, 1],$$
then $\displaystyle D^q_{n,0}(x) =1 $ and for $\displaystyle n > m+ 2$, we have following recurrence relation,
$$[n + m + 2]_qD^q_{n,m+1}(x) = ([m + 1]_q + q^{m+1}x[n]_q)D^q_{n,m}(x) + x(1 - x)q^{m+1}D^q(D^q_{n,m}(x)).$$
\end{theorem}
To establish asymptotic formula for function having $q$-derivative, it is necessary to
 compute moments of first to fourth degree. Using above Theorem one can have first, second, third and fourth order  moments.
\begin{lemma}
For all $x \in [0, 1]$, $n = 1, 2,\ldots$ and $0<q <1$, we have
\begin{itemize}
\item  $\displaystyle D_{n,q}(1,x)=1;$
\item $\displaystyle D_{n,q}(t,x)=\frac{1+q x[n]_q}{[n+2]_q};$
\item $ \displaystyle  D_{n,q}(t^2,x)= \frac{q^3x^2[n]_q[n-1]_q+ (1+q)^2qx[n]_q + 1+ q}{[n+3]_q[n+2]_q};$
\item $ \displaystyle  D_{n,q}(t^3,x)= \frac{q^{8}x^3[n]_q[n-1]_q[n-2]_q  + x^2q^3 [n]_q[n-1]_q\left( 1 + q + 2 q^2 + 3 q^3 + 2 q^4\right)}{[n+4]_q[n+3]_q[n+2]_q}\\
~~~~~~~~~~~~~~~~~~+ \frac{x q[2]_q[n]_q\left( 1 + 2 q + 3 q^2 + 2 q^3 + q^4\right)+ [3]_q [2]_q}{[n+4]_q[n+3]_q[n+2]_q}$;
\item $ \displaystyle D_{n,q}(t^4,x) = \frac{q^{15} x^4[n]_q[n-1]_q[n-2]_q[n-3]_q + q^8 x^3 [n]_q [n-1]_q [n-2]_q \left( 1+2 q+2 q^2+3 q^3+4 q^4+3 q^5+q^6\right)}{[n+5]_q[n+4]_q[n+3]_q[n+2]_q}\\
~~~~~~~~~~~~~~~~~~+ \frac{q^3 x^2 [n]_q[n-1]_q\left\{1+2 q+4 q^2+8 q^3+12 q^4+14 q^5+13 q^6+10 q^7+6 q^8+2 q^9\right\}}{[n+5]_q[n+4]_q[n+3]_q[n+2]_q}\\
~~~~~~~~~~~~~~~~~~+\frac{ q x[2]_q[n]_q\left\{1 + 3 q + 6 q^2 + 9 q^3 + 10 q^4 + 9 q^5 + 6 q^6 + 3 q^7 + q^8\right\} + [4]_q[3]_q [2]_q}{[n+5]_q[n+4]_q[n+3]_q[n+2]_q}$.
 \end{itemize}
\end{lemma}
\begin{lemma}\label{12.lemma1.1}
For all $x \in [0, 1]$, $n = 1, 2,\ldots$ and $0<q <1$, we have
\begin{itemize}
\item $\displaystyle  D_{n,q}\left( (t-x)_q ,x\right) = \frac{1- \left(1+ q^{n+1}\right)x}{[n+2]_q};$
\item $ \displaystyle  D_{n,q}\left( (t-x)_q^2 ,x\right)= \frac{q^2x^2 (1+q^n)(q^{n+1}[2]_q -[n]_q)+x(1+q) (q^2[n]_q   -1-q^{n+2})+ 1+ q}{[n+3]_q[n+2]_q};$
\item $ \displaystyle  D_{n,q}\left( (t-x)_q^3 ,x\right)\\
~~~~~~~~ = q^2 x^3 \left\{\frac{q^{6}[n]_q[n-1]_q[n-2]_q - q[3]_q[n]_q[n-1]_q[n+4]_q + [n+4]_q [n+3]_q[2]_q[n]_q-q[n+4]_q[n+3]_q[n+2]_q}{[n+2]_q[n+3]_q[n+4]_q}\right\}\\
~~~~~~~~+q x^2 \left\{\frac{q^2 [n]_q[n-1]_q\left( 1 + q + 2 q^2 + 3 q^3 + 2 q^4\right) -(1+q)^2[3]_q[n]_q[n+4]_q + [2]_q[n+4]_q [n+3]_q }{[n+2]_q[n+3]_q[n+4]_q}\right\}\\
~~~~~~~~ + x\left\{ \frac{ q[2]_q[n]_q\left( 1 + 2 q + 3 q^2 + 2 q^3 + q^4\right) -(1+ q)[3]_q[n+4]_q }{[n+2]_q[n+3]_q[n+4]_q}\right\} + \frac{[3]_q [2]_q}{[n+2]_q[n+3]_q[n+4]_q};$
\item $ \displaystyle  D_{n,q}\left( (t-x)_q^4 ,x\right)\\
~~~= x^4q^4 \left\{\frac{q^{11}[n]_q[n-1]_q[n-2]_q[n-3]_q}{[n+5]_q[n+4]_q[n+3]_q[n+2]_q} - \frac{q^{4} [4]_q  [n]_q[n-1]_q[n-2]_q}{[n+4]_q[n+3]_q[n+2]_q} + \frac{\left([5]_q+q^2\right)[n]_q[n-1]_q}{[n+3]_q[n+2]_q} - \frac{[4]_q [n]_q}{[n+2]_q}+q^2\right\}\\
~~~~~~+ x^3q^2 \left\{  \frac{q^6  [n]_q [n-1]_q [n-2]_q \left( 1+2 q+2 q^2+3 q^3+4 q^4+3 q^5+q^6\right)}{[n+5]_q[n+4]_q[n+3]_q[n+2]_q} -\frac{ q [4]_q  [n]_q[n-1]_q\left( 1 + q + 2 q^2 + 3 q^3 + 2 q^4\right)}{[n+4]_q[n+3]_q[n+2]_q} \right.\\
~~~~+\left. \frac{(1+q)^2 \left([5]_q+q^2\right)[n]_q}{[n+3]_q[n+2]_q}-q[4]_q \right\}\\
~~~~+ x^2 \left\{ \frac{q^2  [n]_q[n-1]_q\left\{1+2 q+4 q^2+8 q^3+12 q^4+14 q^5+13 q^6+10 q^7+6 q^8+2 q^9\right\} }{[n+5]_q[n+4]_q[n+3]_q[n+2]_q}\right.\\
~~~~\left.-\frac{[4]_q[2]_q [n]_q\left( 1 + 2 q + 3 q^2 + 2 q^3 + q^4\right)}{[n+4]_q[n+3]_q[n+2]_q} +\frac{ (1+ q)\left([5]_q+q^2\right)}{[n+3]_q[n+2]_q}\right\}+\frac{ [4]_q[3]_q [2]_q}{[n+5]_q[n+4]_q[n+3]_q[n+2]_q}\\
~~~~+x\left\{\frac{ q [2]_q[n]_q\left\{1 + 3 q + 6 q^2 + 9 q^3 + 10 q^4 + 9 q^5 + 6 q^6 + 3 q^7 + q^8\right\}+ [4]_q[3]_q [2]_q[n+5]_q}{[n+5]_q[n+4]_q[n+3]_q[n+2]_q}\right\}.$
\end{itemize}
\end{lemma}
\textbf{Proof:} To prove this Lemma, we use linear properties of $q$-Durrmeyer operators.
\begin{eqnarray*}
D_{n,q} ( (t-x)_q , x) &=& D_{n,q}( t,x) - xD_{n,q} (1,x)= \frac{1+q x[n]_q}{[n+2]_q} - x= \frac{1+q x[n]_q- x [n+2]_q}{[n+2]_q} \\
&=& \frac{1+x \left(q+q^2+\ldots+ q^{n}-1-q-q^2-\ldots-q^n-q^{n+1}\right)}{[n+2]_q}\\
& =& \frac{1- \left(1+ q^{n+1}\right)x}{[n+2]_q}.
\end{eqnarray*}
Using identities $\displaystyle (t-x)_q^2 = t^2 - [2]_qxt+qx^2$, we get
\begin{eqnarray*}
D_{n,q} ( (t-x)_q^2 , x) &=& D_{n,q}(t^2,x) -[2]_q xD_{n,q} (t,x)+ qx^2D_{n,q} (1,x)\\
&=&  \frac{q^3x^2[n]_q[n-1]_q+ (1+q)^2qx[n]_q + 1+ q}{[n+3]_q[n+2]_q}
-[2]_qx \left[ \frac{1+q x[n]_q}{[n+2]_q} \right] + qx^2\\
&=&  \frac{q^3x^2[n]_q[n-1]_q+ (1+q)^2qx[n]_q + 1+ q  -[2]_qx [n+3]_q - q x^2[2]_q [n+3]_q [n]_q+  qx^2[n+3]_q[n+2]_q}{[n+3]_q[n+2]_q}\\
&=&  \frac{qx^2 \left\{q^2[n]_q[n-1]_q - [2]_q [n+3]_q [n]_q+ [n+3]_q[n+2]_q\right\}+x\left\{(1+q)^2q[n]_q   -[2]_q [n+3]_q\right\}+ 1+ q}{[n+3]_q[n+2]_q}\\
&=&  \frac{q^2x^2 (1+q^n)(q^{n+1}[2]_q -[n]_q)+x(1+q) (q^2[n]_q   -1-q^{n+2})+ 1+ q}{[n+3]_q[n+2]_q}.
\end{eqnarray*}
Notice that $\displaystyle (t-x)_q^3 = t^3 -[3]_q x t^2 + q[2]_qx^2 t -q^3x^3$,
\begin{eqnarray*}
D_{n,q}\left( (t-x)_q^3 ,x\right)
&=& D_{n,q}\left( t^3 ,x\right)- [3]_q x D_{n,q}\left(  t^2  ,x\right) + q[2]_qx^2 D_{n,q}\left(  t  ,x\right)- q^3x^3\\
&=& \frac{q^{8}x^3[n]_q[n-1]_q[n-2]_q  + x^2q^3 [n]_q[n-1]_q\left( 1 + q + 2 q^2 + 3 q^3 + 2 q^4\right)}{[n+4]_q[n+3]_q[n+2]_q}\\
&&+ \frac{x q[2]_q[n]_q\left( 1 + 2 q + 3 q^2 + 2 q^3 + q^4\right)+ [3]_q [2]_q}{[n+4]_q[n+3]_q[n+2]_q}\\
&&-[3]_q x \left\{ \frac{q^3x^2[n]_q[n-1]_q+ (1+q)^2qx[n]_q + 1+ q}{[n+3]_q[n+2]_q} \right\}+  q[2]_qx^2\left\{ \frac{1+q x[n]_q}{[n+2]_q} \right\}- q^3x^3\\
&&+q x^2 \left\{\frac{q^2 [n]_q[n-1]_q\left( 1 + q + 2 q^2 + 3 q^3 + 2 q^4\right) -(1+q)^2[3]_q[n]_q[n+4]_q + [2]_q[n+4]_q [n+3]_q }{[n+2]_q[n+3]_q[n+4]_q}\right\}\\
&& + x\left\{ \frac{ q[2]_q[n]_q\left( 1 + 2 q + 3 q^2 + 2 q^3 + q^4\right) -(1+ q)[3]_q[n+4]_q }{[n+2]_q[n+3]_q[n+4]_q}\right\} + \frac{[3]_q [2]_q}{[n+2]_q[n+3]_q[n+4]_q}.
\end{eqnarray*}
Finally, using identities $\displaystyle (t-x)_q^4 = t^4-[4]_qx t^3 +q\left([5]_q+q^2\right)x^2 t^2  - q^3x^3 [4]_qt +q^6 x^4$, we get
\begin{eqnarray*}
&&D_{n,q}\left( (t-x)_q^4 ,x\right)\\
&=& D_{n,q}\left( t^4 ,x\right) -[4]_qx D_{n,q}\left( t^3 ,x\right) + q\left([5]_q+q^2\right)x^2  D_{n,q}\left( t^2 ,x\right) -q^3x^3 [4]_q D_{n,q}\left( t ,x\right) +q^6 x^4 \\
&=&\frac{q^{15} x^4[n]_q[n-1]_q[n-2]_q[n-3]_q + q^8 x^3 [n]_q [n-1]_q [n-2]_q \left( 1+2 q+2 q^2+3 q^3+4 q^4+3 q^5+q^6\right)}{[n+5]_q[n+4]_q[n+3]_q[n+2]_q}\\
&&+ \frac{q^3 x^2 [n]_q[n-1]_q\left\{1+2 q+4 q^2+8 q^3+12 q^4+14 q^5+13 q^6+10 q^7+6 q^8+2 q^9\right\}}{[n+5]_q[n+4]_q[n+3]_q[n+2]_q}\\
&&+\frac{ q x[2]_q[n]_q\left\{1 + 3 q + 6 q^2 + 9 q^3 + 10 q^4 + 9 q^5 + 6 q^6 + 3 q^7 + q^8\right\} + [4]_q[3]_q [2]_q}{[n+5]_q[n+4]_q[n+3]_q[n+2]_q}\\
&& - [4]_qx\left\{ \frac{q^{8}x^3[n]_q[n-1]_q[n-2]_q  + x^2q^3 [n]_q[n-1]_q\left( 1 + q + 2 q^2 + 3 q^3 + 2 q^4\right)}{[n+4]_q[n+3]_q[n+2]_q}\right.\\
&&+\left. \frac{x q[2]_q[n]_q\left( 1 + 2 q + 3 q^2 + 2 q^3 + q^4\right)+ [3]_q [2]_q}{[n+4]_q[n+3]_q[n+2]_q}\right\}\\
&&+ q\left([5]_q+q^2\right)x^2 \left\{ \frac{q^3x^2[n]_q[n-1]_q+ (1+q)^2qx[n]_q + 1+ q}{[n+3]_q[n+2]_q} \right\} -q^3x^3 [4]_q\left\{\frac{1+q x[n]_q}{[n+2]_q} \right\}+  q^6 x^4\\
&=&\frac{q^{15} x^4[n]_q[n-1]_q[n-2]_q[n-3]_q + q^8 x^3 [n]_q [n-1]_q [n-2]_q \left( 1+2 q+2 q^2+3 q^3+4 q^4+3 q^5+q^6\right)}{[n+5]_q[n+4]_q[n+3]_q[n+2]_q}\\
&&+ \frac{q^3 x^2 [n]_q[n-1]_q\left\{1+2 q+4 q^2+8 q^3+12 q^4+14 q^5+13 q^6+10 q^7+6 q^8+2 q^9\right\}}{[n+5]_q[n+4]_q[n+3]_q[n+2]_q}\\
&&+\frac{ q x[2]_q[n]_q\left\{1 + 3 q + 6 q^2 + 9 q^3 + 10 q^4 + 9 q^5 + 6 q^6 + 3 q^7 + q^8\right\} + [4]_q[3]_q [2]_q}{[n+5]_q[n+4]_q[n+3]_q[n+2]_q}\\
&& - \frac{q^{8} [4]_q x^4 [n+5]_q [n]_q[n-1]_q[n-2]_q  + q^3 [4]_q x^3 [n+5]_q[n]_q[n-1]_q\left( 1 + q + 2 q^2 + 3 q^3 + 2 q^4\right)}{[n+5]_q[n+4]_q[n+3]_q[n+2]_q}\\
&&- \frac{q[4]_q[2]_qx^2 [n+5]_q[n]_q\left( 1 + 2 q + 3 q^2 + 2 q^3 + q^4\right)+ x [4]_q[3]_q [2]_q[n+5]_q}{[n+5]_q[n+4]_q[n+3]_q[n+2]_q}
\end{eqnarray*}
\begin{eqnarray*}
&&+  \frac{q^4 \left([5]_q+q^2\right)x^4[n+5]_q[n+4]_q [n]_q[n-1]_q+ q^2(1+q)^2 \left([5]_q+q^2\right)x^3[n+5]_q[n+4]_q [n]_q }{[n+5]_q[n+4]_q[n+3]_q[n+2]_q}\\
&&+ \frac{ (1+ q)q\left([5]_q+q^2\right)x^2[n+5]_q[n+4]_q}{[n+5]_q[n+4]_q[n+3]_q[n+2]_q}\\
&& -\frac{q^3x^3[4]_q [n+5]_q[n+4]_q[n+3]_q +q^4x^4 [4]_q [n+5]_q[n+4]_q[n+3]_q[n]_q- q^6 x^4[n+5]_q[n+4]_q[n+3]_q[n+2]_q}{[n+5]_q[n+4]_q[n+3]_q[n+2]_q}\\
&=& x^4q^4 \left\{\frac{q^{11}[n]_q[n-1]_q[n-2]_q[n-3]_q}{[n+5]_q[n+4]_q[n+3]_q[n+2]_q} - \frac{q^{4} [4]_q  [n]_q[n-1]_q[n-2]_q}{[n+4]_q[n+3]_q[n+2]_q} + \frac{\left([5]_q+q^2\right)[n]_q[n-1]_q}{[n+3]_q[n+2]_q} - \frac{[4]_q [n]_q}{[n+2]_q}+q^2\right\}\\
&& + x^3q^2 \left\{  \frac{q^6  [n]_q [n-1]_q [n-2]_q \left( 1+2 q+2 q^2+3 q^3+4 q^4+3 q^5+q^6\right)}{[n+5]_q[n+4]_q[n+3]_q[n+2]_q} -\frac{ q [4]_q  [n]_q[n-1]_q\left( 1 + q + 2 q^2 + 3 q^3 + 2 q^4\right)}{[n+4]_q[n+3]_q[n+2]_q} \right.\\
&&+\left. \frac{(1+q)^2 \left([5]_q+q^2\right)[n]_q}{[n+3]_q[n+2]_q}-q[4]_q \right\}\\
&&+ x^2 \left\{ \frac{q^2  [n]_q[n-1]_q\left\{1+2 q+4 q^2+8 q^3+12 q^4+14 q^5+13 q^6+10 q^7+6 q^8+2 q^9\right\} }{[n+5]_q[n+4]_q[n+3]_q[n+2]_q}\right.\\
&&\left.-\frac{[4]_q[2]_q [n]_q\left( 1 + 2 q + 3 q^2 + 2 q^3 + q^4\right)}{[n+4]_q[n+3]_q[n+2]_q} +\frac{ (1+ q)\left([5]_q+q^2\right)}{[n+3]_q[n+2]_q}\right\}+\frac{ [4]_q[3]_q [2]_q}{[n+5]_q[n+4]_q[n+3]_q[n+2]_q}\\
&&+x\left\{\frac{ q [2]_q[n]_q\left\{1 + 3 q + 6 q^2 + 9 q^3 + 10 q^4 + 9 q^5 + 6 q^6 + 3 q^7 + q^8\right\}+ [4]_q[3]_q [2]_q[n+5]_q}{[n+5]_q[n+4]_q[n+3]_q[n+2]_q}\right\}.
\end{eqnarray*}
\begin{theorem} \label{12.theorem1}Let $f $ bounded and integrable on the interval $[0,1]$ and $(q_n)$ denote a sequence such that $0 < q_n < 1$ and $q_n\to 1$ as $n\to\infty$. Then we have for a point $x \in (0,1)$
$$ \lim_{n\to \infty} [n]_{q_n}[{D}_{n,q_n}(f;x)-f(x)]=(1-2x)\lim_{n\to \infty} D_{q_n}f(x) +x(1-x)\lim_{n\to \infty} D_{q_n}^2f(x).$$
\end{theorem}
\textbf{Proof:  }
By $q$-Taylor formula \cite{de2003integral} for $f$, we have
$$f(t) = f (x)+  D_qf(x) (t - x)+ \frac{1}{[2]_q} D^2_qf(x)(t- x)^2_q+ \theta_q(x;t) (t-x)^2_q,$$
for $0<q<1$, where
\begin{equation}\label{e12}\theta_q(x;t) =\left\{ \begin{array}{cc}
     \displaystyle \frac{f(t) -f(x) -D_qf(x) (t-x) - \frac{1}{[2]_q}D_q^2 f(x) (t-x)^2_q}{(t-x)^2_q} & \text{ if } x\neq t \\
     0, & \text{ if } x=t.
   \end{array}\right.\end{equation}
   We know that for $n$ large enough
\begin{equation}\label{e13}\lim_{t\to x} \theta_q(x; t)=0.\end{equation}
That is for any $\epsilon > 0$, there exists a $\delta > 0$ such that
\begin{equation}\label{e14}|\theta_q(x; t)|\leq \epsilon.\end{equation}
for $|t - x| < \delta$ and $n$ sufficiently large. Using \eqref{e12}, we can write
$${D}_{n,q_n}(f;x)-f(x) = D_{q_n}f(x)  {D}_{n,q_n}((t-x)_q;x) + \frac{D_{q_n}^2f(x)}{[2]_{q_n}}{D}_{n,q_n}((t-x)^2_q;x) + E_n^{q_n}(x),$$
where
$$E_n^{q}(x) =[n+1]_q \sum_{k=0}^n q^{-k} p_{nk}(q; x)\int_0^1 \theta_q(x;t) p_{nk}(q; qt) \left(t-x\right)_q^2 d_qt.$$
By Lemma \ref{12.lemma1.1}, we have
$$ \lim_{n\to \infty} [n]_{q_n} {D}_{n,q_n}((t-x)_q;x)= (1-2x) \text { and } \lim_{n\to \infty} [n]_{q_n} {D}_{n,q_n}((t-x)^2_q;x) = 2 x(1-x).$$
In order to complete the proof of the theorem, it is sufficient to show that $\displaystyle \lim_{n \to \infty}[n]_{q_n}E^{q_n}_n(x) =0$. We proceed as follows:\\
Let
$$P_{n,1}^{q_n}(x) =[n]_{q_n}[n+1]_{q_n} \sum_{k=0}^n {q_n}^{-k} p_{nk}(q_n; x)\int_0^1 \theta_{q_n}(x;t) p_{nk}(q_n; q_nt) \left(t-x\right)_{q_n}^2 \chi_x(t) d_{q_n}t$$
and
$$P_{n,2}^{q_n}(x) =[n]_{q_n}[n+1]_{q_n} \sum_{k=0}^n {q_n}^{-k} p_{nk}(q_n; x)\int_0^1 \theta_{q_n}(x;t) p_{nk}(q_n; q_nt) \left(t-x\right)_{q_n}^2\left(1- \chi_x(t) \right)d_{q_n}t,$$
so that
$$ [n]_{q_n} E_n^{q_n}(x) = P_{n,1}^{q_n}(x) + P_{n,2}^{q_n}(x), $$
where $\chi_x(t)$ is the characteristic function of the interval $\{t:|t - x| < \delta\}$.\\
It follows from \eqref{e12}
$$P_{n,1}^{q_n}(x) = 2 \epsilon  x(1-x) \text { as } n\to \infty.$$
If $\displaystyle |t - x|\geq \delta $, then $\displaystyle |\theta_{q_n}(x; t) |\leq \frac{M}{\delta^2}(t -x)^2$, where $M > 0$ is a constant. Since
\begin{eqnarray*}
\left(t-x\right)^2 &= &\left(t-q^2x + q^2 x -x\right) \left(t-q^3x + q^3 x -x\right)\\
&=& \left(t-q^2x \right) \left(t-q^3x \right)+  x(q^3-1) \left(t-q^2x \right) +x(q^2-1) \left(t-q^2x \right) + x^2(q^2-1)(q^2 - q^3)+  x^2(q^2- 1)(q^3 - 1),
\end{eqnarray*}
we have
\begin{eqnarray*} | P_{n,2}^{q_n}(x) |&\leq& \frac{M}{\delta^2}\left\{ [n]_{q_n}{D}_{n,q_n}((t-x)_{q_n}^4;x) + x(2- q_n^2 -q_n^3)[n]_{q_n} {D}_{n,q_n}((t-x)_{q_n}^3;x) \right. \\
&&\left.+ x^2(q_n^2-1)^2[n]_{q_n} {D}_{n,q_n}^{\alpha,\beta}((t-x)_{q_n}^2;x)\right\}.
\end{eqnarray*}
Using Lemma \ref{12.lemma1.1}, we have
$$ {D}_{n,q_n}((t-x)_{q_n}^4;x) \leq \frac{C_m}{[n]_{q_n}^3},~~~{D}_{n,q_n}((t-x)_{q_n}^3;x) \leq \frac{C_m}{[n]_{q_n}^2} ~~\text{ and } ~~ {D}_{n,q_n}((t-x)_{q_n}^2;x) \leq \frac{C_m}{[n]_{q_n}},$$
we have the desired result.
\begin{corollary} \label{12.coro1}Let $f $ bounded and integrable on the interval $[0,1]$ and $(q_n)$ denote a sequence such that $0 < q_n < 1$ and $q_n\to 1$ as $n\to\infty$. Suppose that the first and second derivative $f'(x)$ and $f''(x)$ exist at a point $x \in (0,1)$. Then we have for a point $x \in (0,1)$
$$ \lim_{n\to \infty} [n]_{q_n}[{D}_{n,q_n}(f;x)-f(x)]=(1-2x)f'(x) +x(1-x) f''(x).$$
\end{corollary}
\section{Asymptotic formula for Durrmeyer-Stancu Operators}
In year 1968, Stancu \cite{stancu1968approximation} generalized Bernstein operators and discussed it approximation properties. After that numbers of researchers gives Stancu type generalization of several operators on finite and infinite intervals, we refer to the papers \cite{mishrasome2013,mishra2013approximation2,mishra2014durrmeyer,buyukyazici2010stancu,verma2012some}. As mention in the introduction Stancu generalization of $q$-Durrmeyer operators \eqref{12.eq1.2} was discussed by Mishra and Patel \cite{mishra2013short}, which is defined as follows: for $0\leq \alpha\leq \beta$,
\begin{equation}\label{12.eq1.3}
D_{n,q}^{\alpha,\beta}=[n+1]_q \sum_{k=0}^n q^{-k} p_{nk}(q; x)\int_0^1 f\left(\frac{[n]_qt+\alpha}{[n]_q+\beta}\right)p_{nk}(q; qt)d_qt,
\end{equation}
where $ \displaystyle p_{nk}(q; x)$ as same as defined in \eqref{12.eq1.2}.
\begin{lemma}\label{l1}
We have
$\displaystyle D_{n,q}^{\alpha,\beta}(1; x)=1,~~~~~\displaystyle D_{n,q}^{\alpha,\beta}(t; x)=  \frac{[n]_q+\alpha[n+2]_q+qx[n]_q^2}{[n+2]_q\left([n]_q+\beta\right)},$\\
 $\displaystyle D_{n,q}^{\alpha,\beta}(t^2; x)= \frac{q^3[n]_q^3\left([n]_q-1\right)x^2+\left(\left(q(1+q)^2+2\alpha q^4\right)[n]_q^3+2\alpha q[3]_q[n]_q^2\right)x}{([n]_q+\beta)^2[n+2]_q[n+3]_q} \\
 \indent \indent \indent \indent \indent +\frac{\alpha^2}{([n]_q+\beta)^2}
 +\frac{(1+q+2\alpha q^3)[n]_q^2+2\alpha[3]_q[n]_q}{([n]_q+\beta)^2[n+2]_q[n+3]_q}.$
 \end{lemma}
 \begin{lemma}\label{l4}
We have \\
$\displaystyle D_{n, q}^{\alpha, \beta}(t-x,x)= \left(\frac{q[n]_q^2}{[n+2]_q([n]_q+\beta)}-1\right)x+\frac{[n]_q+\alpha[n+2]_q}{[n+2]_q([n]_q+\beta)},$\\
$\displaystyle  D_{n, q}^{\alpha, \beta}((t-x)^2,x)= \frac{q^4[n]_q^4-q^3[n]_q^3-2q[n]_q^2[n+3]_q([n]_q+\beta)+[n+2]_q[n+3]_q([n]_q+\beta)^2}{([n]_q+\beta)^2[n+2]_q[n+3]_q}x^2\\
\indent \indent \indent \indent \indent \indent \indent \indent +\frac{q(1+q)^2[n]_q^3+2q\alpha [n]_q^2[n+3]_q-\left(2[n]_q+2\alpha[n+2]_q\right)[n+3]_q([n]_q+\beta)}{([n]_q+\beta)^2[n+2]_q[n+3]_q}x\\
\indent \indent \indent \indent \indent \indent \indent \indent+\frac{(1+q)[n]_q^2+2\alpha[n]_q[n+3]_q}{([n]_q+\beta)^2[n+2]_q[n+3]_q}.$ \end{lemma}
\begin{remark}\label{r4}
For all $m \in\mathbf{N} \cup \{0\} ,0\leq \alpha \leq \beta$; we have the following recursive relation for the images of
the monomials $t^m$ under $D_{n,q}^{\alpha,\beta} (t^m; x)$ in terms of $D_{n,q}(t^j ; x); j = 0, 1, 2,\ldots, m$, as
\begin{eqnarray*}
D_{n,q}^{\alpha,\beta}(t^m;x) =  \sum_{j=0}^{m} {m \choose j}\frac{[n]_q^j\alpha^{m-j}}{([n]_q+\beta)^m}D_{n,q}(t^j,x).
\end{eqnarray*}
 \end{remark}
\begin{theorem}\label{thm2}
Let $f$ bounded and integrable on the interval $[0,1]$ and $(q_n)$ denote a sequence such that $0 < q_n < 1$ and $q_n\to 1$ as $n\to\infty$. Then we have for a point $x \in (0,1)$
$$ \lim_{n\to \infty} [n]_{q_n}[{D}_{n,q_n}^{\alpha,\beta}(f;x)-f(x)]=(1+\alpha-(2+\beta)x)\lim_{n\to \infty} D_{q_n}f(x) +x(1-x)\lim_{n\to \infty} D_{q_n}^2f(x).$$
\end{theorem}
The proof of the above lemma follows along the lines of Theorem \ref{12.theorem1}, using Lemma \ref{l4} and remark \ref{r4}; thus, we omit the details.
\begin{corollary}\label{12.coro2}\cite{mishra2013short}~~Let $f $ bounded and integrable on the interval $[0,1]$ and $(q_n)$ denote a sequence such that $0 < q_n < 1$ and $q_n\to 1$ as $n\to\infty$. Suppose that the first and second derivative $f'(x)$ and $f''(x)$ exist at a point $x \in (0,1)$. Then we have for a point $x \in (0,1)$
$$ \lim_{n\to \infty} [n]_{q_n}[{D}_{n,q_n}^{\alpha,\beta}(f;x)-f(x)]=(1+\alpha-(2+\beta)x)f'(x) +x(1-x)f''(x).$$
\end{corollary}
\begin{remark}
Theorem \ref{12.theorem1} and Theorem \ref{thm2}, gives asymptotic formula for
$q$-Durrmeyer operators and $q$-Durrmeyer-Stancu operators respectively. If $f$ has first and second derivative,
 then $\displaystyle \lim_{n\to \infty} D_{q_n}f(x)=f'(x)$ and $\displaystyle \lim_{n\to \infty} D_{q_n}^2f(x)=f''(x)$.
 We archived results of Mishra and Patel \cite[Theorem 5]{mishra2013short}, which is mention in corollary \ref{12.coro2}.
 So presented results are more general results then exists ones.
\end{remark}


\end{document}